\documentclass[11t]{amsart}
\usepackage{amsaddr}
\usepackage[utf8]{inputenc}
\usepackage{geometry}
\usepackage{ amssymb }
\usepackage[style=numeric,backend=biber]{biblatex}
\usepackage{fancyvrb}
\usepackage{hyperref}
\usepackage{placeins}
\usepackage{tabularx,tabulary}
\usepackage{mathtools}
\usepackage{xcolor}
\definecolor{light-gray}{gray}{0.95}
\newcommand{\code}[1]{\colorbox{white}{\texttt{#1}}}
\usepackage{setspace}
\usepackage{tikz-cd}
\hypersetup{
    colorlinks,
    citecolor=blue,
    filecolor=black,
    linkcolor=blue,
    urlcolor=blue
}

\newtheorem{theorem}{Theorem}

\newtheorem{proposition}[theorem]{Proposition}
\newtheorem{question}[theorem]{Question}
\newtheorem{example}[theorem]{Example}
\newtheorem{claim}[theorem]{Claim}

\theoremstyle{definition}
\newtheorem{definition}[theorem]{Definition}

\theoremstyle{remark}
\newtheorem*{remark}{Remark}

\title[]{Strongly isospectral hyperbolic 3-manifolds with nonisomorphic rational cohomology rings}

\author{ Anda Tenie}
\email{ast2175@columbia.edu}
\address{Department of Mathematics, Columbia University, New York, NY}
\date{\today}

\bibliography{bibliography.bib}

\begin{document}

\begin{abstract}
This paper shows that one cannot ``hear'' the rational cohomology ring of a hyperbolic 3-manifold. More precisely, while it is well-known that strongly isospectral manifolds have the same cohomology as vector spaces, we give an example of compact hyperbolic 3-manifolds that are strongly isospectral but have nonisomorphic rational cohomology rings. Along the way we implement a computer program which finds the nullity of the cup product map $H^1(M;\mathbb{Q})\wedge H^1(M;\mathbb{Q})\rightarrow H^2(M;\mathbb{Q})$ for any aspherical space in terms of the presentation of the fundamental group.

\noindent\emph{\keywordsname:} Isospectral, hyperbolic 3-manifolds, rational cohomology ring, cup product map. 
\end{abstract}

\maketitle 

\section{Introduction}

Two compact Riemannian manifolds $(M,g)$ and $(M',g')$ are said to be \emph{strongly isospectral} if the spectra of every natural self-adjoint elliptic differential operator on $M$
and $M'$ agree (counting multiplicities). If the spectrum of their Laplace operator is the same just on functions we say that the two manifolds are \emph{isospectral}.

Since isometric manifolds are automatically isospectral we are interested in pairs of isospectral manifolds that are not isometric. The construction of these has a rich history
starting with the example of Milnor \cite{milnor1964eigenvalues} in 1964 of 16-dimensional isospectral flat tori. 
While this example cannot be generalized in some natural way, in 1985 Sunada found a general method of constructing strongly isospectral manifolds, also known as the Sunada method \cite{sunada1985riemannian}. This construction also led to the famous example by Gordon, Webb, and Wolper \cite{gordon1992one} of regions in the plane that are isospectral but not isometric, hence, answering Kac's question ``Can you hear the shape of a drum?" negatively.

We are interested in better understanding how much topology two strongly isospectral hyperbolic 3-manifolds have to share. It is known that in general two strongly isospectral manifolds have the same Betti numbers and, therefore,
$$H^*(M_1; \mathbb{Q})=H^*(M_2; \mathbb{Q})$$ as vector spaces (this follows from the Hodge decomposition theorem). The following question is then natural:

\begin{question}
Do strongly isospectral manifolds have isomorphic rational cohomology rings? In other words, does their ring structure have to be the same?
\end{question} 
In 2013, Lauret, Miatello, and Rossetti \cite{lauret2013strongly} answered the question negatively by exhibiting pairs of compact high-dimensional flat manifolds that are strongly isospectral, having nonisomorphic cohomology rings. On the other hand, we are concerned with hyperbolic 3-manifolds. These spaces are particularly interesting because of their role in low-dimensional topology and number theory. The main goal of this paper is to present an example of strongly isospectral hyperbolic 3-manifolds with nonisomorphic cohomology rings. 

\begin{theorem}
There exists a pair of strongly isospectral compact hyperbolic 3-manifolds with nonisomorphic rational cohomology rings.
\end{theorem}

These examples are obtained via the Sunada method which we review in Section \ref{sunadamethod}. We then find a pair with distinct cohomology rings  by studying the maps $$H^1(M;\mathbb{Q})\wedge H^1(M;\mathbb{Q})\rightarrow H^2(M;\mathbb{Q}).$$
The strategy is to find examples for which these maps have kernels of different ranks. In Section \ref{complexker}, we show it suffices to compute the kernels of $$\mu_K: H^1(K;\mathbb{Q})\wedge H^1(K;\mathbb{Q})\rightarrow H^2(K;\mathbb{Q})$$ where $K$ is the 2-complex associated to some finite presentation of the fundamental groups of our manifolds. In order to compute this map, we follow the paper of Suciu and Wang \cite{suciu2019cup} which, for a 2-complex $K,$ gives a description of the cup-product map $H^1(K;\mathbb{Q})\wedge H^1(K;\mathbb{Q})\rightarrow H^2(K;\mathbb{Q}).$ We explain this method, as well as its implementation, in Section \ref{cupproductsection}. Along the way, we give a computer program that finds the above cup product map for any aspherical space. Finally, in Section \ref{exampleiso}, we present the example of two strongly isospectral compact hyperbolic 3-manifolds with nonisomorphic rational cohomology rings.

\section{The Sunada Method and its implementation}
\label{sunadamethod}

We begin by briefly describing the Sunada method. For a more detailed treatment which includes the proof of the Sunada theorem one can consult Sunada's original paper \cite{sunada1985riemannian} or Buser's text \cite[Chapter 11]{buser2010geometry}. We then explain the implementation of this construction which, given a hyperbolic manifold from the SnapPy census, produces a pair of strongly isospectral manifolds as covering spaces of the input (if such a pair exists).

\begin{definition}
 Let $G$ be a finite group and $H_1$ and $H_2$ two subgroups of $G.$ We say $H_1$ and $H_2$ are \emph{almost conjugate} if for all $g,$
$$|[g]\cap H_1|=|[g]\cap H_2|.$$ In other words, the two subgroups intersect each conjugacy class in the same number of elements.
\end{definition}

\begin{example}\label{example}
Some examples of finite groups $G$ with $H_1$ and $H_2$ almost conjugate are:

\begin{enumerate}
\item $G=PSL(2,7)$ with two nonconjugate index 7 subgroups.
\item $G=\mathbb{Z}_8^{*}\ltimes \mathbb{Z}_8$ with multiplication $(a,b)(a',b')=(aa',ab'+b)$ and  $$H_1=\{(1,0);(3,0);(5,0);(7,0)\}$$ 
$$H_2=\{(1,0);(3,4);(5,4);(7,0)\}$$
\end{enumerate}
\end{example}
\begin{theorem}[Sunada]  Let $M$ be a complete Riemannian manifold and $G$ a finite group acting on $M$ by isometries with at most finitely many fixed points. If $H_1$ and $H_2$ almost conjugate acting freely on $M$ then the quotients $M/H_1$ and $M/H_2$ are strongly isospectral.
\end{theorem}

\begin{remark}Sunada's theorem does not claim that the two quotients are not isometric. In order to make sure the two manifolds are not isometric, we looked for strongly isospectral manifolds with different integral homology. Such examples with different torsion in integral homology had already been shown to exist (Bartel, Page \cite{bartel2016torsion}).\end{remark}

Sunada's theorem then gives a general method for finding pairs of strongly isospectral manifolds. The algorithm is as follows:

Consider a group $G$ with almost conjugate subgroups $H_1$ and $H_2$ and a hyperbolic 3-manifold $M_0$ whose fundamental group surjects into $G$: $$\varphi:\pi_1(M_0) \twoheadrightarrow{} G\supset H_1,H_2.$$ Let $M$ be the covering space corresponding to the kernel of this map and $M_1, M_2$ the covering spaces corresponding to the preimages $\varphi^{-1}(H_1) $ and $\varphi^{-1}(H_2).$ If we let $M_1=M/H_1$ and $M_2= M/H_2$ then by Sunada's theorem $M_1$ and $M_2$ are strongly isospectral.
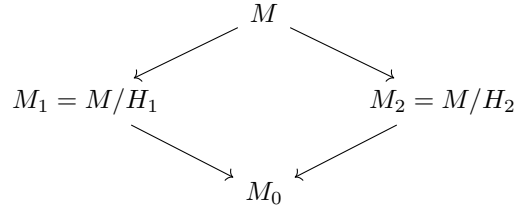
\[\begin{tikzcd}
	{} & M \\
	{M_1=M/H_1} && {M_2=M/H_2} \\
	& {M_0}
	\arrow[from=1-2, to=2-1]
	\arrow[from=2-1, to=3-2]
	\arrow[from=1-2, to=2-3]
	\arrow[from=2-3, to=3-2]
\end{tikzcd}\]

\begin{remark}
First examples of strongly isospectral hyperbolic 3-manifolds were constructed by Reid \cite{reid1992isospectrality} in 1992  using Sunada's method.
\end{remark}

In this work, we primarily used the group $G=\mathbb{Z}_8^{*}\ltimes \mathbb{Z}_8$ with almost conjugate subgroups $H_1, H_2$ as in Example \ref{example}.
In GAP, the group $G$ can be identified with SmallGroup(32,43). If we let \code{C:=ConjugagyClassesSubgroups(G)} then representatives of \code{C[14]} and \code{C[15]} are almost conjugate. We show this by checking that every element in \code{C[14]} is conjugate to an element of \code{C[15]} as done below. This code can also be found in the file AlmostConjugateExample.ipynb available \href{https://github.com/andatenie/CupAspherical}{here}.

\begin{verbatim}
G=gap.SmallGroup(32,43)
G=gap.SimplifiedFpGroup(gap.Image(gap.IsomorphismFpGroup(G)))
C=gap.ConjugacyClassesSubgroups(G)
gap.RelatorsOfFpGroup(G)
Output: [ F1^2, F2^2, F3^2, (F3*F2)^2, (F1*F2*F1*F3)^2, F2*F1*(F2*F1*F3)^2*F2*F1 ]
C[15].Representative()
Output: Group( [ F2, F1*F3*F1^-1 ] )
C[14].Representative()
Output: Group( [ F2, F3 ] )
gap.ConjugacyClass( G, G.2*G.1*G.3*G.1)==gap.ConjugacyClass( G, G.2*G.3)
Output: True
G.1*G.2*G.1*(G.2*G.1*G.3*G.1)*G.1*G.2*G.1==G.2*G.3
Output: True
\end{verbatim}

As seen above, the representative of \code{C[15]} is generated by $f_2$ and $f_1f_3f_1$ and the representative of \code{C[14]} is generated by $f_2$ and $f_3.$ So, to show that the two groups are almost conjugate, it suffices to show that $f_2 \cdot f_1f_3f_1$ is in the same conjugacy class as $f_2\cdot f_3.$ This is precisely what is done in the last two lines of the above code.

\begin{remark}
These are the same almost conjugate subgroups as the ones presented in Example \ref{example} with identifications: $f_1=(7,1), f_2=(7,0),$ and $f_3=(5,0).$
\end{remark}

With the almost conjugate subgroups constructed as above we can now look for corresponding strongly isospectral covers. Given a hyperbolic 3-manifold $M_0$ from SnapPy, our code checks if there is a surjection $\varphi:\pi_1(M_0) \twoheadrightarrow{}G.$ If such a surjection is found, a pair of covering spaces $M_1$ and $M_2$ is constructed and, by Sunada's theorem, the two manifolds are strongly isospectral. Below is the code that was used for this construction which can also be found in the file SunadaPair.ipynb available \href{https://github.com/andatenie/CupAspherical}{here}. In order to run it one needs SageMath to work together with SnapPy. The tool we used is \href{https://github.com/sagemath/docker-images}{SageMath Docker Image}. 
\begin{verbatim}
S=gap.SmallGroup(32,43);
C=gap.function_call('ConjugacyClassesSubgroups',S);
H1=gap.function_call('Representative',C[15]);
H2=gap.function_call('Representative',C[14]);
M=snappy.Manifold("t10238(0,0)")
M.dehn_fill((2,1),0);
G=M.fundamental_group();
G=gap(G);
if(G.GQuotients(S).Length()>0):
    f = G.GQuotients(S)[1];
    H1p=gap.PreImage(f,H1);
    H2p=gap.PreImage(f,H2);
    M1 = M.cover(H1p);
    M2 = M.cover(H2p);
\end{verbatim}
The base manifold chosen here, \code{t10238(2,1)} is, as we will see in Section \ref{exampleiso}, the one which provides the pair $M_1, M_2$ with nonisomorphic rational cohomology rings. The rest of the paper is devoted to showing that $H^*(M_1)\ncong H^*(M_2).$
\section{The cup product map}\label{complexker}
We now turn our attention to understanding the cup product map $$\mu: H^1(M;\mathbb{Q})\wedge H^1(M;\mathbb{Q})\rightarrow H^2(M;\mathbb{Q}).$$
We want to find examples of a Sunada pair of hyperbolic 3-manifolds $M_1, M_2$ whose kernels have different ranks and, hence, different rational cohomology rings.

 Let $K$ be the 2-complex associated to some finite presentation of $\pi(M),$ where the generators are the 1-cells and the relators are the corresponding 2-cells.
Consider the cup product map $$\mu_K: H^1(K;\mathbb{Q})\wedge H^1(K;\mathbb{Q})\rightarrow H^2(K;\mathbb{Q}).$$

\begin{claim}
The maps $\mu$ and $\mu_K$ have the same kernel for $M$ aspherical. Hence, it suffices to understand the latter.
\end{claim}

\begin{remark}
The importance of the above claim comes from the fact that $\mu_K$ can be computed using the method of Suciu and Wang \cite{suciu2019cup}.
\end{remark}

\begin{proof}
 Note that we have an inclusion $K \hookrightarrow K'=  K(\pi_1(M),1)$ into an Eilenberg–MacLane space $K'$ formed by attaching cells of dimension $\geq 3$ on $K$ to kill higher homotopy groups. Therefore, we obtain a diagram

\[\begin{tikzcd}
	{H^1(K')\wedge H^1(K')} & {H^2(K')} \\
	{H^1(K)\wedge H^1(K)} & {H^2(K)}
	\arrow["\mu", from=1-1, to=1-2]
	\arrow["f"',  from=1-1, to=2-1]
	\arrow["{\mu_K}"', from=2-1, to=2-2]
	\arrow["g", from=1-2, to=2-2]
\end{tikzcd}\]
Since $M= \mathbb{H}^3/\pi_1(M),$ we note that $M$ is also a $K(\pi_1(M),1)$ space. Since a weak homotopy equivalence induces an isomorphism on cohomology we get 

\[\begin{tikzcd}
	{H^1(M)\wedge H^1(M)} & {H^2(M)} \\
	{H^1(K)\wedge H^1(K)} & {H^2(K)}
	\arrow["\mu", from=1-1, to=1-2]
	\arrow["f"',  from=1-1, to=2-1]
	\arrow["{\mu_K}"', from=2-1, to=2-2]
	\arrow["g", from=1-2, to=2-2]
\end{tikzcd}\]
Note that $f$ is an isomorphism since $\pi_1(M)=\pi_1(K).$ Moreover, $g$ is an injection since $H^2(K',K)=0.$ Therefore, $\text{ker}(\mu)=\text{ker}(\mu_K).$ 
\end{proof}
\section{Cup product map for 2-complexes}\label{cupproductsection}
 Given a presentation of $\pi_1(M)$ one can construct a 2-complex $K$ corresponding to that presentation. In \cite{suciu2019cup}, Suciu and Wang give a formula for the map $\mu_K: H^1(K;\mathbb{Q})\wedge H^1(K;\mathbb{Q})\rightarrow H^2(K;\mathbb{Q}).$ We recall this method for the reader's convenience and then explain how we implemented the calculation of this map.
 
 \subsection{Fox Derivatives and their implementation}

Let $F$ be a free group on generators ${x_1,\dots, x_n}$ and $\mathbb{Z}F$ its group ring.
The \emph{Fox derivatives} $\partial_i:\mathbb{Z}F \rightarrow \mathbb{Z}F$ satisfy the following rules:  for all $v_1, v_2 \in \mathbb{Z}F,$ $$\partial_i(v_1 +v_2)=\partial_i(v_1 )+\partial_i(v_2)$$ $$\partial_i(v_1 v_2)=\partial_i(v_1 )\varepsilon(v_2)+v_1\partial_i(v_2),$$
where $\varepsilon$ is the augmentation map $\varepsilon(\sum g_in_i)=\sum n_i,$ for $ g_i \in F.$ Note that if $g_1$ and $g_2$ are elements of the group then $$\partial_i(g_1g_2)=\partial_i(g_1)+g_1\partial_i(g_2).$$

The following are useful consequences of this definition: $\partial_i(1)=0,~ \partial_i(x_j)=\delta_{ij},$ and $\partial_i(g^{-1})=-g^{-1}\partial_i(g)$ for all $g$ in $F.$ The higher Fox derivatives $\partial_{i_1,\dots, i_k}$ are defined inductively. Let the operator $\varepsilon_I: \mathbb{Z}G\rightarrow \mathbb{Z}$ be the composition $\varepsilon_I=\varepsilon \circ \partial_I.$ For more details on Fox derivatives the reader can consult \cite{fox1953free},\cite[Chapter VII]{crowell1964introduction}.

We implemented the Fox derivative as a recursive function called \code{FoxDerivative}. It takes in \code{element} (the word to be differentiated), \code{index} (the index with respect to which we take the derivative), and \code{result} (a variable which keeps track of the result so far in the recursion). Since the group presentations use letters of the alphabet, we encode generators of the group as $\{a,b,c,\dots\}$ and their respective inverses as $\{A,B,C,\dots\}.$ The index will then be a letter in $\{a,b,c, \dots\}.$ The function will read the first letter of \code{element} and compare it to \code{index} in order to understand the case we are in. For example, $\partial_a(aBc)=1+a\partial_a(Bc)$ in which case the function is called again on \code{element}$=Bc$ and \code{index}$=a,$ the result of this being multiplied by $a$ after which we add 1. On the other hand note, $\partial_a(BBc)=B\partial_a(Bc)$ is a different case. The function \code{FoxDerivative} handles all such cases and returns the final result as an element of the free algebra on the generators in the presentation and their inverses.
\begin{remark}Being a recursion, this function becomes less efficient as the words of the presentation get longer. For this reason, if one only wants to compute the composition of the augmented map with the Fox derivative i.e. $\varepsilon\circ\partial_i,$ the function \code{AugmentFoxDerivative} can be used. Note that this loses information about the first derivative so it cannot be used to determine higher Fox derivatives. In order to compute the double Fox derivatives (which will be needed in determining the cup product map) we will compute the first using \code{FoxDerivative} and the second one using \code{AugmentFoxDerivative} in order to improve efficiency. 
\end{remark}
\code{AugmentFoxDerivative} takes in \code{element}(the word to be differentiated), \code{index}(the index with respect to which we take the derivative). It is based on the following observation: $\varepsilon\circ\partial_\text{\code{index}}(\text{\code{element}})$ counts the number of times \code{index} appears in the word \code{element}, adding 1 every time it sees \code{index} and subtracting one every time it sees its inverse. This happens because $$\partial_a (al_2l_3\dots l_k)=1+a \partial_a (l_2l_3\dots l_k) $$  $$\partial_a (Al_2l_3\dots l_k)= -A + A\partial_a (l_2l_3\dots l_k).$$

For example, $\varepsilon\circ \partial_a(abAcAB)=-1$ since it contains one copy of $a$ and two of $A.$ This can also be seen by direct computation:
$\varepsilon(1 - abA - abAcA)= 1-1-1=-1.$
\subsection{Echelon presentations}
\begin{definition}
 Let $G$ be a finite group with finite presentation $ \langle \mathbf{x}~ |~ \mathbf{ w} \rangle.$
 We say this is an \emph{echelon presentation} if the augmented Fox Jacobian matrix $(\varepsilon_i(w_k))$ is in row echelon form. 
\end{definition}
The above definition is significant because we will be able to compute the cup products for echelon presentations. The following proposition explains how to construct an echelon presentation of any finitely presented group $G.$

\begin{proposition} \label{echelonexists}
Let $G$ a group with finite presentation $ \langle x_1,\dots x_n |~ r_1,\dots r_m \rangle.$ There exists then a group $G_e$ with echelon
presentation, and a map $f : K_{G_e} \rightarrow K_G$ between the respective presentation 2-complexes such that the induced
homomorphism in cohomology, $f^*:H^*(K_G;\mathbb{Z})\rightarrow H^*(K_{G_e};\mathbb{Z}),$ is an isomorphism.
\end{proposition}

For a proof of this proposition see \cite[Proposition 3.3]{suciu2019cup}. The group $G_e$ is constructed as follows: consider the transpose of the $n\times m$ matrix $\varepsilon_i(r_k)$ and denote this by $T.$ Then there exists some $C=(c_{l,k})$ invertible such that $C \cdot T$ is in Hermite normal form. The group $G_e$ with echelon presentation is given by $$G_e=\langle x_1,\dots x_n |~ w_1,\dots w_m \rangle$$ where $w_k=r_1^{c_{1,k}}r_2^{c_{2,k}}\dots r_m^{c_{m,k}}$ for $1 \leq k \leq m. $ 

We used the \code{FoxDerivative} function described above to construct the transpose of $\varepsilon_i(r_k).$ GAP's function \code{HermiteNormalFormIntegerMatTransform} can then find the matrix $C$ such that $C \cdot T$ is in Hermite normal form. Upon finding this matrix it is simple to construct the new relators as in the formula of $w_k$ given above. The new relators are kept in a list called \code{RelatorsEchelon}. The code can be found in the file CupProduct.ipynb available \href{https://github.com/andatenie/CupAspherical}{here}.

\subsection{Computing the cup product map}

Let $G$ be a group with finite presentation $\langle x_1,\dots,x_n ~| ~r_1,\dots r_m \rangle$ and let $F\twoheadrightarrow{} G$ be the map from the free group on  $x_1,\dots,x_n $ to $G.$ Then, as we have seen in {Proposition~\ref{echelonexists}}, there exists $G_e$ with echelon presentation given by $\langle x_1,\dots,x_n ~ | ~w_1,\dots,w_m  \rangle$ and map $f:K_{G_e}\rightarrow K_G$ inducing an isomorphism in cohomology. Let $\{u_1,\dots, u_b\}$ and $\{\beta_{n-b+1},\dots, \beta_{m}\}$ be the basis for $H^1(K;\mathbb{Q})$ and $H^2(K;\mathbb{Q})$ respectively, obtained by  transfering bases in $H^*(K_G;\mathbb{Q})$ from suitable bases for $H^*(K_{G_e};\mathbb{Q}).$ The following gives an explicit formula for the cup product map.

\begin{theorem}[\cite{suciu2019cup}]\label{thmSW}
Then the cup product map $\mu_K: H^1(K;\mathbb{Q})\wedge H^1(K;\mathbb{Q})\rightarrow H^2(K;\mathbb{Q})$ is given by  $$u_i\cup u_j=\sum_{k=n-b+1}^m \kappa(w_k)_{i,j}\beta_k$$
where $\kappa$ is the Magnus $\kappa$-extension for $F$ relative to $G.$\end{theorem}
Moreover, the coefficients $\kappa(w_k)_{i,j}$ can be computed using the following result:
\begin{proposition}[\cite{suciu2019cup}]
Let $(a_{i,s})$ be the $b\times n$ matrix associated to the linear map $\pi: F_\text{ab}\otimes \mathbb{Q}\rightarrow G_\text{ab}\otimes \mathbb{Q}.$ Then for each $1\leq i,j\leq b$ we get $$\kappa(w_k)_{i,j}=\sum_{s,t=1}^n a_{i,s}a_{j,t}\varepsilon_{s,t}(r).$$
\end{proposition}

We begin by explaining how the program computes the matrix $(a_{i,s}).$ Note that $G_{ab}$ is the cokernel of $\mathbb{Z}\langle w_1,\dots,w_m\rangle\rightarrow \mathbb{Z}\langle x_1,\dots,x_n\rangle $ where the map takes a word to the linear combination of the letters that are part of it. For example, it would take $x_1x_1x_2x_2x_2+x_4$ to $2x_1+3x_2+x_4.$ We can then compute the matrix corresponding to this map using the \code{AugmentFoxDerivative} function described before and then find the $(a_{i,s})$ as the cokernel.

Once $(a_{i,s})$ are found the only pieces needed are $\varepsilon_{s,t}(r).$ These are computed using the function \\\code{DoubleFoxDeriv} which takes the double Fox derivative by first applying \code{FoxDerivative} and then\\ \code{AugmentFoxDeriv}. The final step consists in building the matrix $(\kappa(w_k)_{i,j})$ for $1\leq i<j\leq b$ and $n-b+1\leq k \leq m$ and finding the kernel of the map. The details of this construction can be found in the file CupProduct.ipynb available \href{https://github.com/andatenie/CupAspherical}{here}.
\section{The example}\label{exampleiso}
\begin{theorem}
There exists a pair of strongly isospectral compact hyperbolic 3-manifolds with nonisomorphic rational cohomology rings.
\end{theorem}

The pair of isospectral manifolds is constructed using the Sunada method recalled in Section \ref{sunadamethod} using the group $G=\mathbb{Z}_8^{*}\ltimes \mathbb{Z}_8$ with $H_1$ and $H_2$ almost conjugate given as in Example \ref{example}. The base manifold was found after an extensive search and is obtained by Dehn filling (2,1) the manifold \code{t10238(0,0)}, which is number 15053 in the  OrientableCuspedCensus. The two strongly isospectral manifolds $M_1,M_2$ that we obtain have the following integral homology:

$$H_1(M_1; \mathbb{Z})=\mathbb{Z}/2 \oplus\mathbb{Z}/2 \oplus\mathbb{Z}/2 \oplus\mathbb{Z}/2 \oplus\mathbb{Z}/4 \oplus \mathbb{Z}/4  \oplus \mathbb{Z} \oplus \mathbb{Z} \oplus \mathbb{Z}$$
$$H_1(M_2; \mathbb{Z})=\mathbb{Z}/2 \oplus\mathbb{Z}/2 \oplus\mathbb{Z}/2 \oplus\mathbb{Z}/2 \oplus\mathbb{Z}/2 \oplus \mathbb{Z}/8  \oplus \mathbb{Z} \oplus \mathbb{Z} \oplus \mathbb{Z}.$$

Hence, they are not isometric. In order to compute their cup product map we first need their fundamental groups:

$\pi_1(M_1)=\langle$\code{a,b,c,d,e,f,g,h,i| ahAIcGBGHcHicahbAc, aHACbfhgadifcHcbHdAgCi, acgbgCAB,}\newline
\code{aeBeAb, aefdAgCicaheIDAGHF, ahAcahAfcHcbHFBc, ahBHAgbg, adhChCdAgChChg, adhbHdAB}$\rangle$

$\pi_1(M_2)=\langle$\code{a,b,c,d,e,f,g,h,i| ccdechCiFeiFeHde, bEifbCiFeC, aafDifaafggDifaaff, }\newline \code{aiaffaafgDDifaafgf,bHbchC, bCEDbcdecciFif, aafbcdecbEdGGf, bCiFeiFifbCEfI, aiAi}$\rangle.$

Now let $K_1$ and $K_2$ be the 2-complexes associated to these presentations. Note that dim$H^1(K_i;\mathbb{Q})=b=3$ and so dim$H^1(K_i;\mathbb{Q})\wedge H^1(K_i;\mathbb{Q})=3$ for $i\in\{1,2\}.$ Moreover, by Theorem \ref{thmSW}, dim$H^2(K_i;\mathbb{Q})=3$ for $i\in\{1,2\}.$ The cup product map $H^1(K_1;\mathbb{Q})\wedge H^1(K_1;\mathbb{Q})\rightarrow H^2(K_1;\mathbb{Q})$ is in some basis given by
$$\begin{pmatrix}
-2 & -2 & 2 \\
0 & 2 & -2 \\
2 & 2 & 0
\end{pmatrix}$$ and has kernel of rank 0. On the other hand, $H^1(K_2;\mathbb{Q})\wedge H^1(K_2;\mathbb{Q})\rightarrow H^2(K_2;\mathbb{Q})$ is given by 
$$\begin{pmatrix}
0 & 0 & 0 \\
0 & 0 & 0 \\
0 & 0 & 0
\end{pmatrix}$$
 and so it has kernel of rank 3. This shows that the rank of the kernels is different and so the maps $$\mu: H^1(M_i;\mathbb{Q})\wedge H^1(M_i;\mathbb{Q})\rightarrow H^2(M_i;\mathbb{Q})$$ also have kernels of different rank. Hence, the two strongly isospectral manifolds have nonisomorphic rational cohomology rings.
 
 \begin{remark}
  A theorem of Thurston \cite{gromov1981hyperbolic} says that the values of the function $V\mapsto\text{Vol}(V)$ form a closed nondiscrete well-ordered set on the real line where $V$ runs over all hyperbolic 3-manifolds of finite volume. It is then natural to ask: what is the smallest volume that an isospectral pair can have? Previous work of Linowitz and Voight \cite{linowitz2015small} using number theoretic methods found a volume $\approx51.02$ example. The pair we give above has the smallest volume found so far, $\approx 25.418347$, and we believe this would be difficult to improve significantly due to certain constraints. For example, the index of the covering spaces is least 7 and the smallest volume of a closed hyperbolic 3-manifold is $\approx 0.9427,$ corresponding to the Weeks manifold \cite{gabai2009minimum}.
 \end{remark}

\section*{Acknowledgement}
  I would like to thank my advisor, Professor Francesco Lin, for his guidance throughout this research project, for introducing me to this research question, and for sharing with me his knowledge and enthusiasm for the field. This research was supported by the Alfred P. Sloan Foundation.
\newpage
\nocite{*}
\printbibliography

\end{document}